\documentclass[11pt]{article}
\input epsf
\usepackage{amsfonts}
\usepackage{amscd}
\usepackage{amsmath,amssymb,amsthm}
\usepackage{amstext}
\usepackage{amscd}
\usepackage{graphicx}
\usepackage[all]{xy}
\usepackage{pstricks,pst-node,pst-tree}
\usepackage[T1]{fontenc}   
\usepackage{lmodern}

\theoremstyle{plain}
\newtheorem{lem}{Lemma}[section]
\newtheorem*{cubical}{Theorem A}
\newtheorem*{mfld}{Theorem B}
\newtheorem*{main1}{Theorem 1}
\newtheorem*{main2}{Corollary 2}

\theoremstyle{definition}

\newtheorem{definition}[lem]{Definition}

\newtheorem{rem}[lem]{Remark}

\renewcommand{\descriptionlabel}[1]%
       {\hspace{\labelsep}\textsf{#1}}

\newcommand{\R}{\mathbb{R}}

\begin{document}
\title{A Constructive Cubical Realization of $n$-Dimensional Smooth Knots Inside the Menger 
$M^{n+2}_n$-continuum\,\,}
\author{Juan Pablo D\'iaz, Gabriela Hinojosa,  Alberto Verjovsky\thanks{This work was partially supported by PAPIIT (Universidad 
Nacional Aut\'onoma de M\'exico) project \#IN103324.}}
\date{May 15, 2026}

\maketitle

\begin{abstract} 
We prove that every smooth $n$-dimensional knot in $\mathbb{R}^{n+2}$ can be ambiently isotoped into the Menger $n$-dimensional continuum. In contrast with classical embedding theorems for universal compacta, our construction is explicit and proceeds via cubical models, combining the cubical realization theorem of Boege--Hinojosa--Verjovsky with the affine self-similarity of the Menger continuum.
\end{abstract}

\section{Introduction}

It is classical that the  Menger $M^{n+2}_n$-continuum is universal for all compact metric spaces (compacta) of topological dimension smaller or equal to $n$ that embed in $\mathbb{R}^{n+2}$, for instance by results of \v{S}tan'ko 
\cite{Stanko} and Bestvina \cite{Bestvina1988}.  Such results are existential in nature.\\

\noindent The purpose of this paper is different.  We give a constructive cubical
realization of smooth codimension-two knots inside the standard Menger
continuum.  The construction proceeds in three steps:
\begin{enumerate}
    \item isotoping the knot into the canonical cubical scaffolding ($n$-skeleton),
    \item realizing this scaffolding inside a cubical Menger model $M_N$,
    \item transporting $M_N$ ambiently onto a finite iteration of the  Menger compactum $M^{n+2}_n$.
\end{enumerate}

\noindent Thus, the embedding is obtained explicitly through cubical isotopies and affine
self-similarity.\\

\noindent More precisely,  we will prove the following.

\begin{main1}
Let $N^n$ be a cubical closed $n$-dimensional submanifold embedded in the $n$-skeleton of the canonical cubulation ${\mathcal{C}}^{n+2}$ 
of $\mathbb{R}^{n+2}$. Then there exists an isotopic copy of  $N^n$ into the  Menger $M^{n+2}_n$-continuum. 
\end{main1}

\noindent Boege--Hinojosa--Verjovsky \cite{BHV} proved that any smooth $n$-dimensional knot is isotopic to  a knot
contained in the codimension-two skeleton (scaffolding) of the canonical cubulation $\mathcal{C}$ of
$\mathbb{R}^{n+2}$ (see section 2). As a consequence, we have the following.

\begin{main2}
Any smooth knot $\mathbb{S}^n\hookrightarrow\mathbb{R}^{n+2}$ is isotopic to an $n$-knot contained in the Menger 
$M^{n+2}_n$-continum.
\end{main2}

\noindent This is a generalization of the theorem proved by  Broden--Espinosa--Nazareth--Voth, which states that all one-dimensional knots in $\mathbb{R}^3$ can also be found within a finite iteration of the Menger sponge (see \cite{Barber}, \cite{BENV}). We would like to emphasize that, if we take $n=1$ in Corollary 2, we obtain the above result; however, our result is not a straight forward generalization since we use completely different arguments.

\section{Preliminaries}

This section introduces the basic concepts and some results about cubical $n$-manifolds and  Menger $M^{n+2}_n$-continuum.\\

\noindent The canonical cubulation $\mathcal{C}$ of $\mathbb{R}^{n+2}$ is the decomposition into hypercubes, which are the images of the unit cube 
$I^{n+2} = \{(x_1,\ldots,x_{n+2}) \mid 0 \le x_i \le 1\}$ by translations by vectors with integer coefficients.\\

\noindent Consider the homothetic transformation $h_m : \mathbb{R}^{n+2} \to \mathbb{R}^{n+2}$ given by $h_m(x) = m x$, where $m > 1$ is an integer. The set $h_m(\mathcal{C})$ is called a subcubulation or cubical subdivision of $\mathcal{C}$.\\

\begin{definition} 
The $k$\emph{-skeleton} of the canonical cubulation ${\mathcal{C}}^{n+2}$ of $\mathbb{R}^{n+2}$, denoted by ${\mathcal{C}}^{n+2}_k$,
consists of the union of the $k$-skeletons of the hypercubes in  ${\mathcal{C}}^{n+2}$,
{\it i.e.,} the union of all cubes of dimension $k$ contained in  the $(n+2)$-cubes in  ${\mathcal{C}}^{n+2}$.
We will call the $n$-skeleton ${\mathcal{C}}^{n+2}_n$ of ${\mathcal{C}}^{n+2}$, denoted by ${\mathcal{S}}^{n}$, the {\it canonical scaffolding} of $\mathbb{R}^{n+2}$.
\end{definition}

\subsection{Cubical submanifolds}

\noindent We will start with some basic facts about cubical submanifolds; for more details, see \cite{DHVV} and \cite{DHV1}.

\begin{definition}
Let $N^n$ be an $n$-dimensional closed submanifold in $\mathbb{R}^{n+2}$. If $N^n$ is contained in the canonical scaffolding $\mathcal{S}^n$, we say that $N^{n}$ is a \emph{cubical $n$-submanifold}. 
\end{definition}

\noindent By \cite{BHV}, we have many examples of cubical $n$-submanifolds. More specifically, we have the following (\cite{BHV}).

\begin{cubical}
Let  $M^{n+1},\,\,N^n\subset\mathbb{R}^{n+2}$, $N^n\subset{M^{n+1}}$,
be closed and smooth submanifolds of  $\mathbb{R}^{n+2}$.
Suppose that $N^n$ has a trivial normal bundle in $M^{n+1}$ ({\it i.e.,} $N^n$ is a two-sided hypersurface of $M^{n+1}$). Then there exists an ambient isotopy of
$\mathbb{R}^{n+2}$ which takes $M^{n+1}$ into the $(n+1)$-skeleton of the canonical
cubulation $\mathcal{C}^{n+2}$ of $\mathbb{R}^{n+2}$ and $N^n$ into the $n$-skeleton of $\mathcal{C}^{n+2}$. In particular, $N^n$ can be deformed by
an ambient isotopy into a cubical manifold contained in the canonical scaffolding of $\mathbb{R}^{n+2}$.
\end{cubical}

\noindent Notice that a cubical $n$-submanifold in $\mathbb{R}^{n+2}$ is a 
piecewise linear submanifold such that each linear piece is a unit $n$-hypercube with its vertices in the $\mathbb{Z}^{n+2}$-lattice of $\mathbb{R}^{n+2}$. 
We remark that, in general, our cubical $n$-submanifolds are, in a natural way, discrete hyperbolic length spaces (\cite{BB}).\\

\noindent We recall that, in classical knot theory, a subset $K$ of a space $X$ is a {\it knot} if $K$ is homeomorphic to a  $p$-dimensional sphere 
$\mathbb{S}^{p}$ embedded in either the Euclidean $n$-space $\mathbb{R}^{n}$ or the $n$-sphere $\mathbb{S}^{n}=\mathbb{R}^{n}\cup\{\infty\}$, where $p<n$. Two knots $K_1$, $K_2$ are 
{\it equivalent} or {\it isotopic} if there is a homeomorphism $h:X\hookrightarrow X$ such that $h(K_1)=K_2$;
in other words $(X,K_1)\cong (X,K_2)$. However, a knot $K$ is sometimes defined to be an embedding
$K:\mathbb{S}^{p}\hookrightarrow\mathbb{S}^{n}\cong\mathbb{R}^{n} \cup \{\infty \}$ (see \cite{mazur}, \cite{rolfsen}).
We shall also find this convenient at times and will use the same symbol to
denote either the map $K$ or its image $K(\mathbb{S}^{p})$ in $\mathbb{S}^{n}$.

\begin{definition}
Let $K^n$ be a $n$-dimensional knot in $\mathbb{R}^{n+2}$. If $K^n$ is contained in the canonical scaffolding $\mathcal{S}^n$, we say that $K^{n}$ is a \emph{cubical $n$-knot}. 
\end{definition}

\noindent In particular, we have the following.

\begin{mfld}[Boege--Hinojosa--Verjovsky \cite{BHV}]
Let $K\subset\mathbb{R}^{n+2}$
be a smooth $n$-dimensional knot.  Then $K$ is ambiently isotopic to a knot
contained in the codimension-two skeleton of the canonical cubulation of
$\mathbb{R}^{n+2}$.  Equivalently, since the ambient dimension is $n+2$, it is
ambiently isotopic to a knot contained in the canonical cubical $n$-skeleton.
\end{mfld}

\begin{rem} This theorem, in dimension 3 has been used extensively and apparently has been proven but we were not able to find a specific citation. Although it is to
be expected in high dimensions, there was no proof in the literature. The proof is technical
and relies on cubulations of tubular neighborhoods together with carefully constructed ambient
isotopies adapted to the cubical structure. It is not a standard PL or cubical approximation
theorem.
\end{rem}

\subsection{The cubical Menger compactum}

Let
\[
Q=[0,1]^{n+2}.
\]
Subdivide $Q$ into $3^{n+2}$ congruent subcubes indexed by
\[
(a_1,\ldots,a_{n+2})\in\{0,1,2\}^{n+2},
\]
where $a_i=0,1,2$ correspond respectively to the intervals
\[
[0,1/3],\qquad [1/3,2/3],\qquad [2/3,1].
\]
The first-stage approximation $M_n^{n+2}(1)$ of the standard $n$-dimensional
Menger continuum is obtained by keeping precisely those closed subcubes for
which at most $n$ of the indices $a_i$ are equal to $1$.  Equivalently, one
removes the open subcubes whose addresses have at least $n+1$ middle
coordinates.\\

\noindent The same rule is then iterated inside every retained cube.  This gives a
nested sequence of compact cubical sets
\[
Q=M_n^{n+2}(1)\supset M_n^{n+2}(1)\supset M_n^{n+2}(2)\supset\cdots,
\]
and the standard Menger continuum is
\[
M^{n+2}_n=\bigcap_{r=0}^{\infty}M_n^{n+2}(r).
\]

\begin{figure}[h] 
 \begin{center}
 \includegraphics[height=4.7cm]{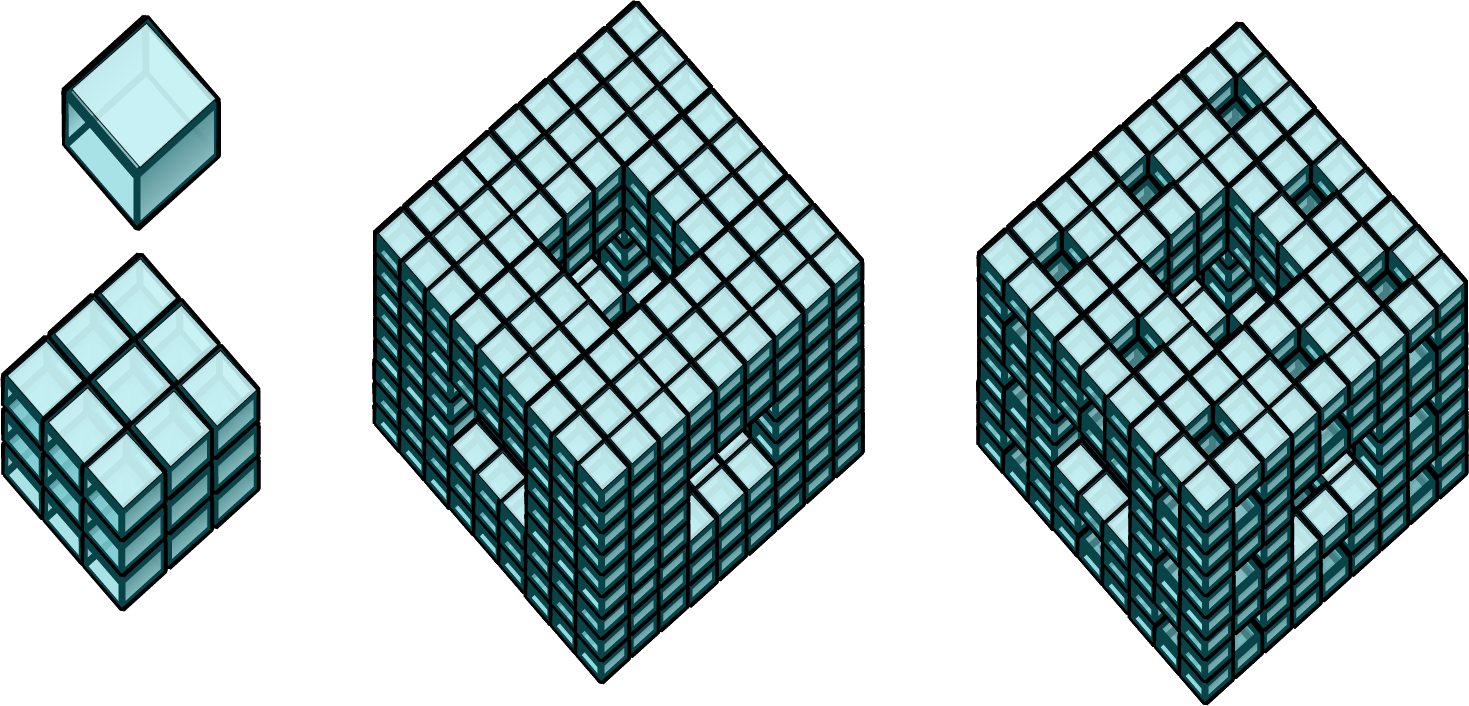}
\end{center}
\caption{\sl First steps in the construction of Menger's sponge.} 
\label{Menger}
\end{figure}

\noindent For $n=1$, this is the classical Menger sponge in $[0,1]^3$: one removes the
central cube and the six face-center cubes at the first stage, and then
iterates in every retained cube (see Figure \ref{Menger}).\\

\noindent The main technical point is the following elementary but crucial observation:
The standard Menger rule removes only cells that are too central in too many
coordinate directions.  Consequently, no cubical cell of dimension at most $n$
is deleted from the retained scaffolding.

\section{Main Theorem}

The goal of this section is to prove Theorem 1.

\begin{main1}
Let $N^n$ be a cubical closed $n$-dimensional submanifold embedded in the $n$-skeleton of the canonical cubulation ${\mathcal{C}}^{n+2}$ 
of $\mathbb{R}^{n+2}$. Then there exists an isotopic copy of  $N^n$ into the  Menger $M^{n+2}_n$-continuum. 
\end{main1}

\begin{figure}[h] 
 \begin{center}
 \includegraphics[height=4.5cm]{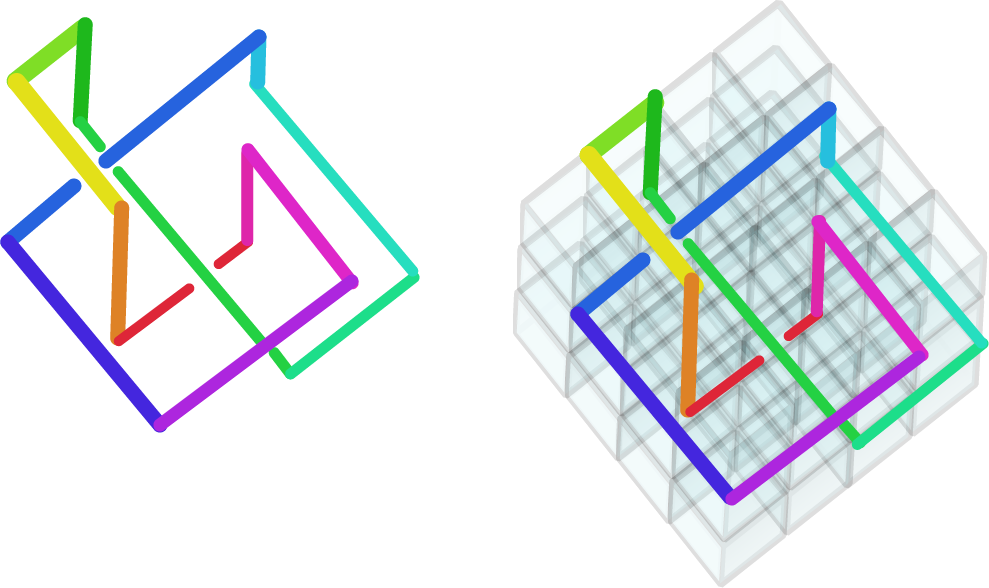}
\end{center}
\caption{\sl $N^n$ embedded in the $n$-skeleton of the hypercube $C^{n+2} (2^{m}-1)$.} 
\label{Paso0}
\end{figure}

\noindent{\it Proof.} First, we will find a proper $(n+2)$-dimensional cube such that our knot $N^n$ is contained in the $n$-skeleton of it (see Figure \ref{Paso0}). Then, we apply a construction to get a space homeomorphic, via a homothetic transformation, to some finite stage  $\tilde{M}^{n+2}_n (k)$ of the Menger $\tilde{M}^{n+2}_n$- continuum,  in such a way that an isotopic copy of $N^n$ is contained in the $n$-skeleton of it. We will start describing our recursive construction.\\

\noindent Since $N^n$ is compact, there is an integer  $m$ such that $N^n$ is embedded in the  $n$-skeleton of an $(n+2)$-dimensional hypercube 
$C^{n+2} (2^{m}-1)$ contained into $\mathbb{R}^{n+2}$ with length side $2^{m}-1$ (see Figure \ref{Paso0}). \\

\begin{figure}[h] 
 \begin{center}
 \includegraphics[height=5cm]{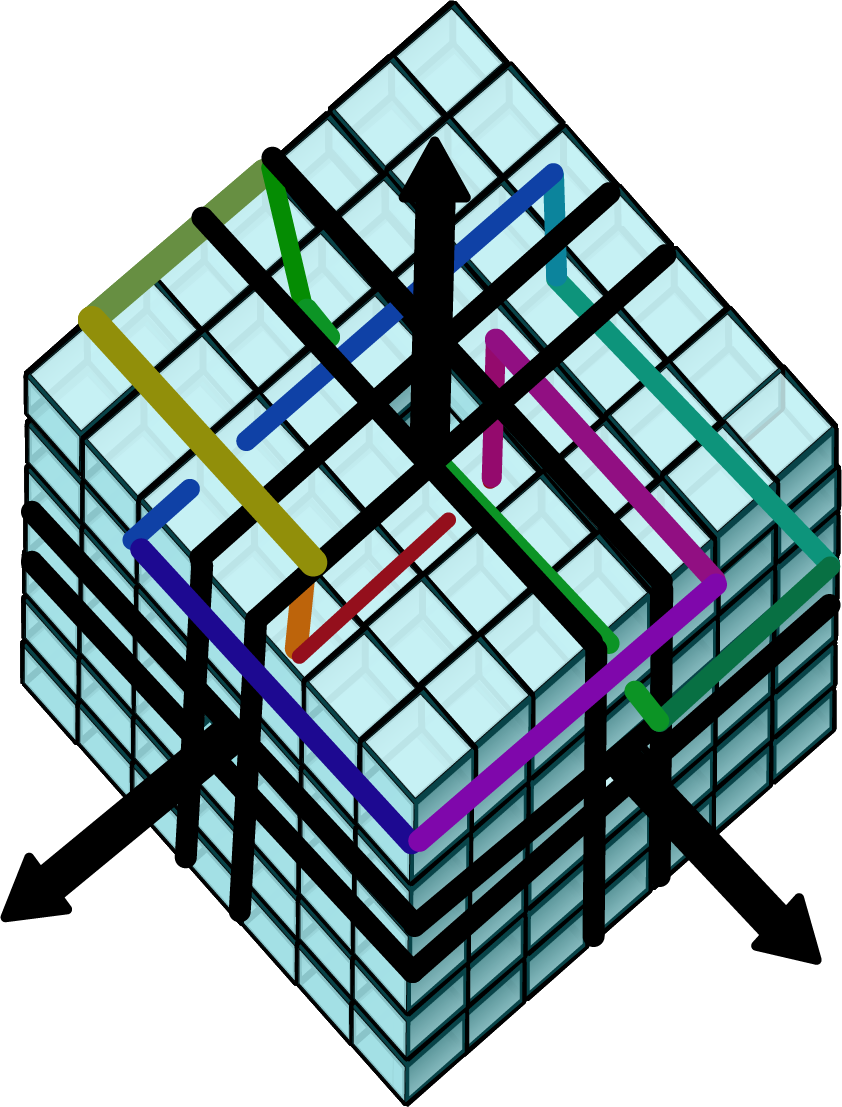}
\end{center}
\caption{\sl The central hypercubes of length side one are removed from  $C^{n+2} (2^{m}-1)$.} 
\label{Paso2}
\end{figure}

\noindent {\it First stage}. We will start with the $(n+2)$-dimensional hypercube $C^{n+2} (2^{m}-1)$. As in the first stage of the construction of the usual $(n+2)$-Menger continuum, we remove the corresponding hypercubes at the center of each $(n+1)$-face and the central hypercube of $C^{n+2} (2^{m}-1)$, but in our case, their side length is one. Then,  we obtain a set $D_1$ that can be decomposed into $2^{n+2}$ hypercubes $H^1_i$ ($i=1,\ldots, 2^{n+2}$) of side length  $2^{m-1}-1$ (one hypercube for each vertex of $C^{n+2} (2^{m}-1)$), and $4\times 3^{n}$ hyperslices (hyper parallelepipeds) $S^1_j$ whose size is
$(2^{m-1}-1) \times \dots \times (2^{m-1}-1) \times 1= \prod_{i=0}^{n-1}  (2^{m-1}-1) \times 1 $, that shares a $(n+1)$-face $F^1_{j,s}$ of length side $2^{m-1}-1$ with  two different hypercubes $H^1_s$, for $j=1,\ldots 4\times 3^{n}$; in other words, any side of these hyperslices has length $2^{m-1}-1$, except one edge whose length is one (see Figure \ref{Paso2}).\\

\noindent We would like to point out that $N^n$ is contained in the $n$-skeleton of $D_1$, since we have only removed the interiors of $(n+2)$-dimensional unit hypercubes of $C^{n+2}_{2^{m}-1}$. \\

\begin{figure}[h] 
 \begin{center}
 \includegraphics[height=4cm]{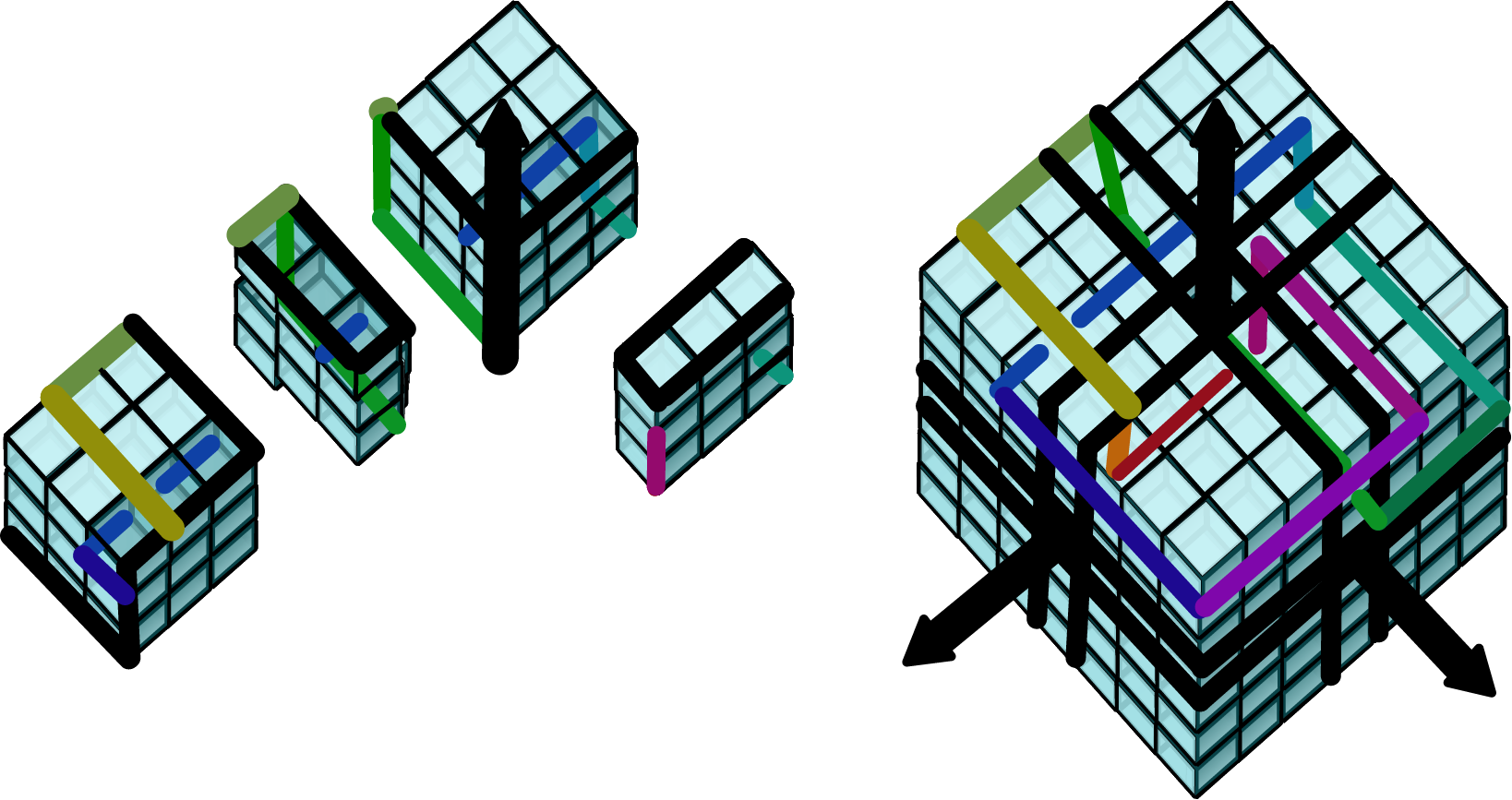}
 \includegraphics[height=4cm]{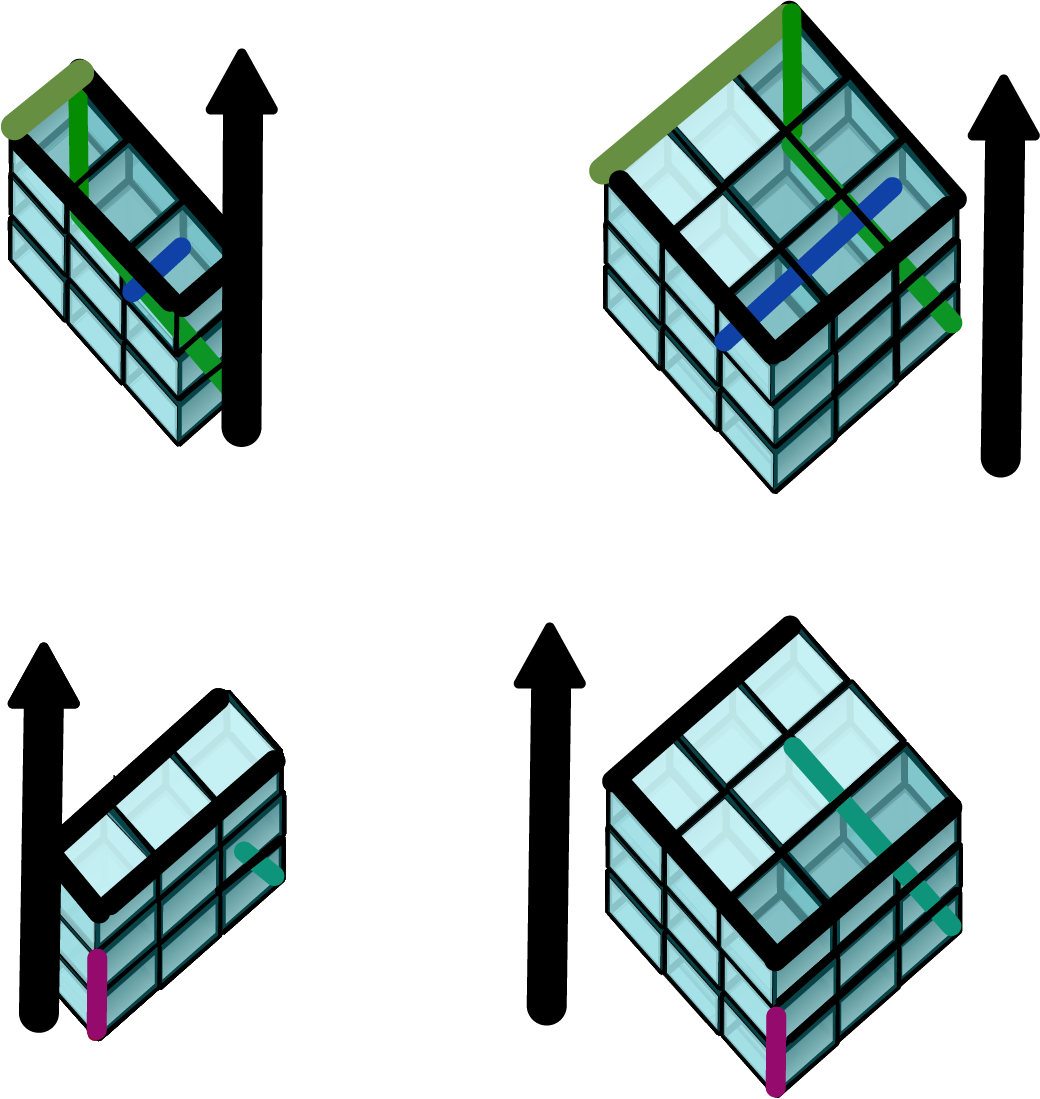}
\end{center}
\caption{\sl The hypercubes $H^1_i$ and the hyperslices (hyper parallelepipeds) $S^1_j$.} 
\label{Paso2ab}
\end{figure}

\noindent Now, we will modify  $D_1$  to obtain a space $B_1$ equivalent to the first step  $M^{n+2}_n (1)$ of the usual construction of the $(n+2)$-Menger continuum, in such a way that an isotopic copy of $N^n$ is contained in the $n$-skeleton of $B_1$ (see Figure \ref{Paso2ab}).
 Indeed, to obtain a scale model of $M^{n+2}_n (1)$, we must rescale the boundary of the removed hypercubes in such a way that their side length will be $2^{m-1}-1$, instead of one. To do this, we multiply the edge of length one on each hyperslice $S^1_j$ by $2^{m-1}-1$, for $j=1,\ldots 4\times 3^{n}$; in other words, we apply a homothetic transformation $T$ to expand this one side length, in such a way that we obtain $\tilde{S}^1_j$ which is a hypercube of side length $2^{m-1}-1$, and the hypercubes $H^1_j$ stay in the same way. Notice that the $(n+1)$-face $F^{1}_{j,s}=S^1_j\cap H^1_r$ is equal to the $(n+1)$-face $\tilde{S}^1_j\cap H^1_s$ for each $j$ and $s$, hence our space $B_1$ is obtained from $D_1$ replacing $S^1_j$ by $\tilde{S}^1_j$ and identifying its corresponding $(n+1)$-faces with the corresponding $(n+1)$-faces of the hypercubes $H^1_s$. Now, the same expansion applies to each intersection $N^n\cap S^1_j$ in such a way that we get an isotopic copy $N^n_1$ of $N^n$ that is contained in the $n$-skeleton of  $B_1$.\\
 
\begin{figure}[h] 
 \begin{center}
 \includegraphics[height=5cm]{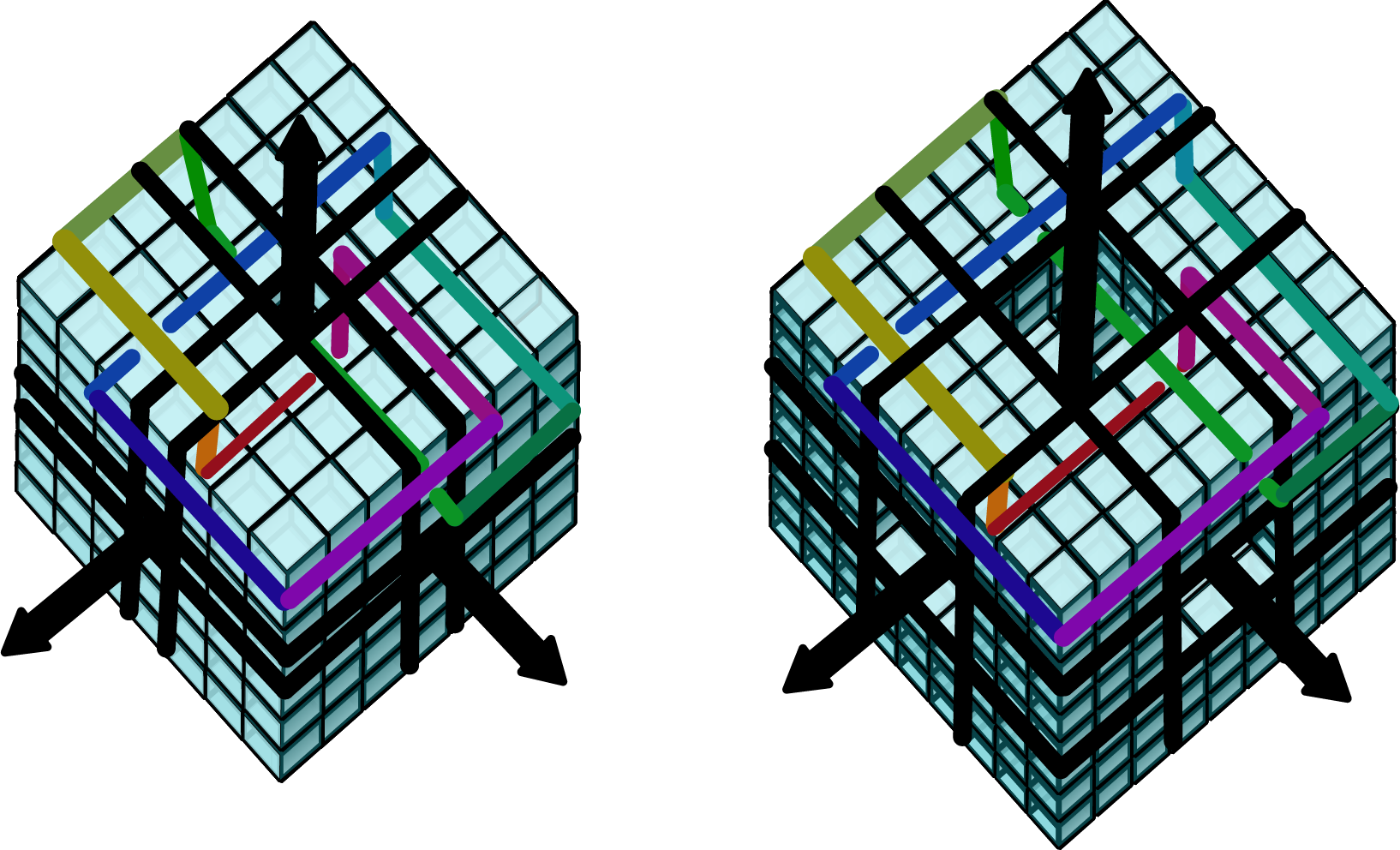}
\end{center}
\caption{\sl First stage of our construction.} 
\label{Paso2c}
\end{figure}

\noindent Observe that by construction, $B_1$ consists of all those subcubes that intersect the $n$-skeleton of our modified hypercube $C^{n+2} (3\times (2^{m-1}-1))$ of side length $3\times (2^{m-1}-1)$, that has been divided into $3^{n+2}$ subcubes with side length $2^{m-1}-1$; in other words, $B_1$ is homothetic to the corresponding space $M^{n+2}_n (1)$ obtained in the first stage of the construction of the $(n+2)$-dimensional Menger $M^{n+2}_n$-continuum (see Figure \ref{Paso2c}).\\

\noindent {\it Second stage.} We continue as in the construction of the Menger $M^{n+2}_n$-continuum. This is, consider each
$(n+2)$-hypercube of length side $2^{m-1}-1$ of $B_1$ and repeat the first stage on it. We again remove the corresponding hypercubes at the center of each $(n+1)$-face and the central hypercube of each hypercube of $B_1$. Again, in our case,  their length side is one, so we get a set $D_2$ consisting of hypercubes $H^2_i$ of length side $2^{m-2}-1$ and hyperlices $S^2_j$ whose length side is $(2^{m-2}-1) \times \dots \times (2^{m-2}-1)  \times 1= \prod_{i=0}^{n-1}  (2^{m-2}-1) \times 1$ that shares a $(n+1)$-face $F^2_{j,s}$ of length side $2^{m-2}-1$ with two different hypercubes $H^2_s$. Again, we expand each $S^2_j$ to get a hypercube 
$\tilde{S}^2_j$ of length side $2^{m-2}-1$, in such a way that we get a space $B_2$ obtained from $D_2$ replacing each $S^2_j$ for the corresponding hypercube 
$\tilde{S}^2_j$ and identifying its corresponding $(n+1)$-faces with the corresponding $(n+1)$-faces of the hypercubes $H^2_s$. As before, the same expansion applies to each intersection $N^n_1\cap S^1_j$, so we get an isotopic copy $N^n_2$ of $N^n_1$ that is contained in the $n$-skeleton of $B_2$.\\

\begin{figure}[h] 
 \begin{center}
 \includegraphics[height=5cm]{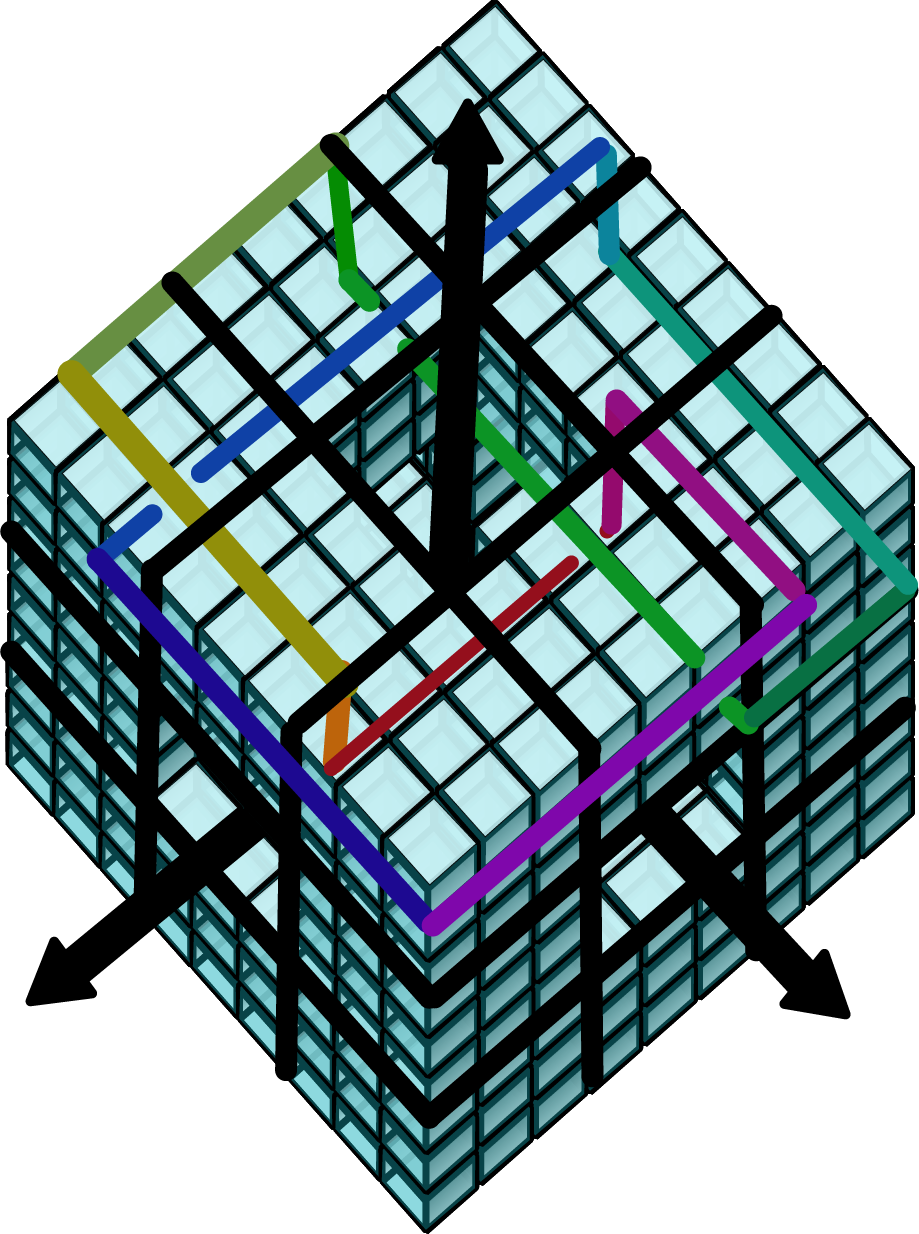}\hspace{2cm}
 \includegraphics[height=5cm]{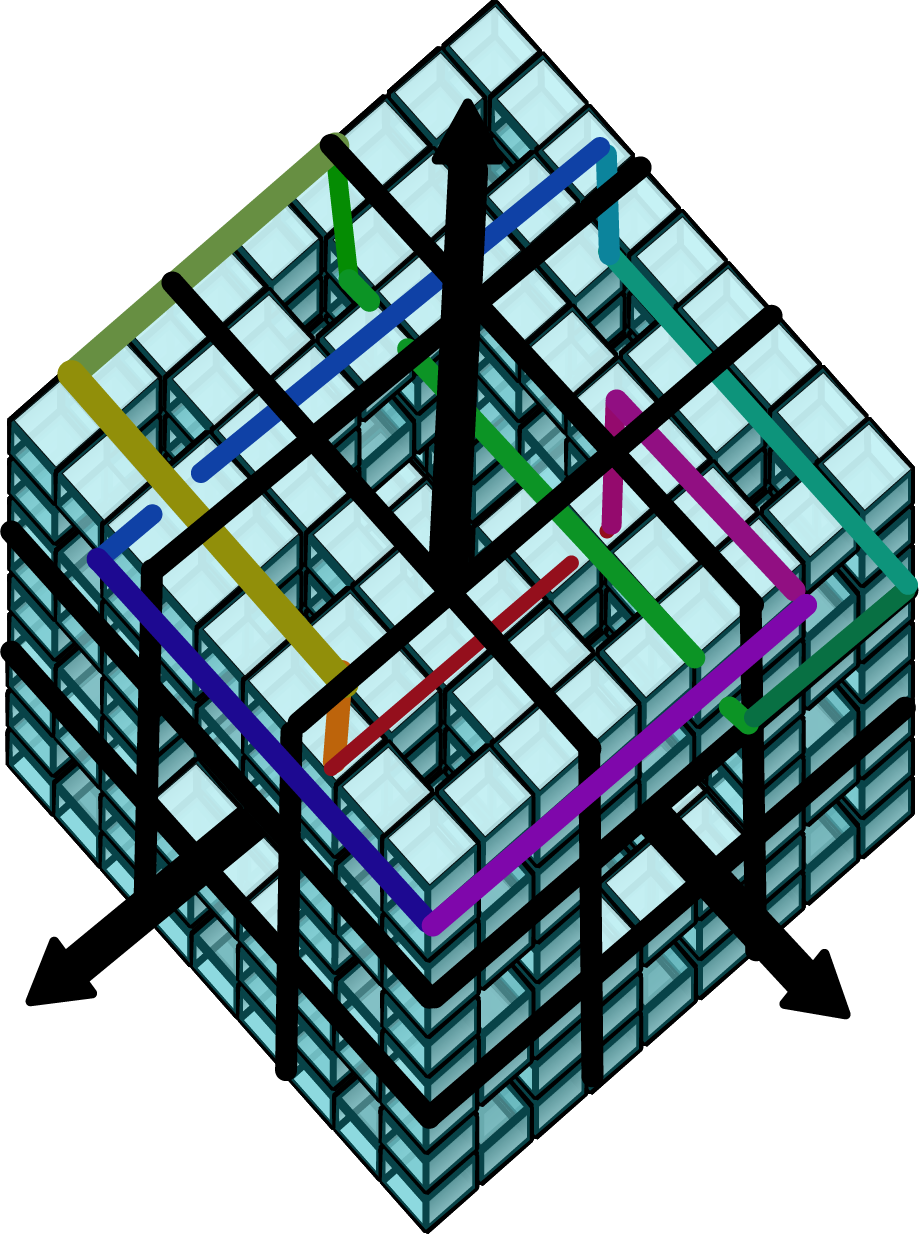}
\end{center}
\caption{\sl Second stage of our construction.} 
\label{Paso3}
\end{figure}

\noindent By construction, $B_2$ consists of all those subcubes that intersect the $n$-skeleton of  hypercubes of $B_1$ whose side length is $2^{m-1}-1$. Notice that each subcube of $B_2$ has side length $2^{m-2}-1$ and the lenght side of our modified hypercube $C^{n+2} (3^2\times (2^{m-2}-1))$ is $3^2\times (2^{m-2}-1)$. Hence $B_2$ is obtained as the second stage $M^{n+2}_n (2)$ of the  Menger $M^{n+2}_n$-continuum construction, hence they are homothetic (see Figure \ref{Paso3}).\\

\noindent {\it $(m-1)$-stage.} We continue as in the construction of the  $M^{n+2}_n$-continuum. This is, consider each
$(n+2)$-hypercube of length side $2^{m-(m-2)}-1$ of $B_{m-2}$ and repeat the first stage on it. We again remove the corresponding hypercubes at the center of each $(n+1)$-face and the central hypercube of each hypercube of $B_{m-2}$ whose their length side is one, to get a set $D_{m-1}$ consisting only of hypercubes $H^{m-1}_i$ of length side $2^{m-(m-1)}-1=1$; in fact, there are no hyperslices at this stage, hence $D_{m-1}$ is the space $B_{m-2}$. As a consequence, the submanifold $N^n_{m-2}$ obtained in the previous stage which  is isotopic to  $N^n$, is contained in the $n$-skeleton of $B_{m-1}$.\\

\noindent Notice that by construction, $B_{m-1}$ consists of all those subcubes that intersect the $n$-skeleton of the hypercube of $B_{m-2}$ whose side length is
 $2^{m-(m-2)}-1$. Notice that each subcube of $B_{m-1}$ has side length $1$ and the lenght side of our modified hypercube $C^{n+2} (3^{m-1})$ is 
 $3^{m-1}\times (2^{m-(m-1)}-1)=3^{m-1}$. Hence $B_{m-1}$ is obtained as the corresponding space $M^{n+2}_n (m-1)$ in the $(m-1)$-stage of the  Menger 
 $M^{n+2}_n$-continuum, hence they are homeomorphic. Even more, we can rescale $B_{m-1}$ just applying the homothetic transformation 
 ${\mathcal{T}}:\mathbb{R}^{n+2}\rightarrow \mathbb{R}^{n+2}$ given by  ${\mathcal{T}}(x) =\frac{1}{3^{m-1}} (x)$, thus 
 ${\mathcal{T}}(C^{n+2} (3^{m-1}) =\frac{1}{3^{m-1}} (C^{n+2} (3^{m-1}) =I^{n+2}$, and ${\mathcal{T}}(B_{m-1}) =\frac{1}{3^{m-1}} (B_{m-1})=M^{n+2}_n (m-1)$.\\


\noindent We continue this process as the construction of the $M^{n+2}_n$-continuum. Notice that an isotopic copy of our cubical, closed $n$-submanifold 
$N^n$ is contained in the Menger $M^{n+2}_n$-continuum, since an isotopic copy of it is contained in $B_{m-1}$. Therefore, the result follows. $\square$

\begin{main2}
Any smooth knot $\mathbb{S}^n\hookrightarrow\mathbb{R}^{n+2}$ is isotopic to an $n$-knot contained in the Menger 
$M^{n+2}_n$-continum.
\end{main2}

\noindent{\it Proof.} Let $\tilde{K}\subset \R^{n+2}$ be a smooth  $n$-dimensional knot in $\mathbb{R}^{n+2}$, then by \cite{BHV}, there exists a cubical $n$-knot $K$ in 
$\mathbb{R}^{n+2}$ isotopic to $\tilde{K}$. The result follows the previous theorem. $\square$

\noindent J. P.  D\'iaz. {\tt Centro de Investigaci\'on en Ciencias}. Instituto de Investigaci\'on en Ciencias B\'asicas y Aplicadas. Universidad Aut\'onoma del Estado de Morelos. Av. Universidad 1001, Col. Chamilpa.
Cuernavaca, Morelos, M\'exico, 62209. 

\noindent {\it E-mail address:} juanpablo.diaz@uaem.mx
\vskip .3cm

\noindent G. Hinojosa. {\tt Centro de Investigaci\'on en Ciencias}. Instituto de Investigaci\'on en Ciencias B\'asicas y Aplicadas. 
Universidad Aut\'onoma del Estado de Morelos. Av. Universidad 1001, Col. Chamilpa.
Cuernavaca, Morelos, M\'exico, 62209. 

\noindent {\it E-mail address:} gabriela@uaem.mx

\vskip .3cm

\noindent A. Verjovsky. {\tt Instituto de Matem\'aticas Unidad Cuernavaca}. Universidad Nacional Aut\'onoma de M\'exico. Av. Universidad 1001, Col. Chamilpa.
Cuernavaca, Morelos, M\'exico, 62210. 

\noindent {\it E-mail address:} albertoverjovsky@gmail.com

\end{document}